\documentclass{article}

\usepackage{arxiv}

\usepackage[utf8]{inputenc} 
\usepackage[T1]{fontenc}    
\usepackage{hyperref}       
\usepackage{url}            
\usepackage{booktabs}       
\usepackage{amsfonts}       
\usepackage{nicefrac}       
\usepackage{microtype}      
\usepackage{lipsum}		
\usepackage{graphicx}
\usepackage{natbib}
\usepackage{doi}

\usepackage{amsmath,amssymb,enumerate,xspace}

\usepackage{booktabs, siunitx}
\usepackage{booktabs, tabularx, amsmath}

\usepackage{amsthm}
\theoremstyle{definition} 

\usepackage{color}
\definecolor{dukeblue}{rgb}{0.0, 0.0, 0.61}
\definecolor{harvardcrimson}{rgb}{0.79, 0.0, 0.09}
\definecolor{dartmouthgreen}{rgb}{0.05, 0.5, 0.06}
\definecolor{gold(metallic)}{rgb}{0.83, 0.69, 0.22}
\definecolor{goldenrod}{rgb}{0.85, 0.65, 0.13}
\definecolor{indianyellow}{rgb}{0.89, 0.66, 0.34}
\definecolor{indianred}{rgb}{0.8, 0.36, 0.36}
\definecolor{eggplant}{rgb}{0.38, 0.25, 0.32}

\title{The Strong Birthday Problem Revisited}

\author{{Chijul B.~Tripathy} \\
	Homestead High School\\
	Cupertino, CA 95014 \\
	\texttt{chijul.b.tripathy@gmail.com} \\
}

\date{}

\renewcommand{\shorttitle}{}

\hypersetup{
pdftitle={The Strong Birthday Problem Revisited},
pdfauthor={Chijul B.~Tripathy},
pdfkeywords={The birthday problem, Strong birthday problem, Stirling number of the second kind, Associated Stirling number of the second kind, Dynamic programming, Recurrence},
}

\newcommand{\sbp}{\texttt{SBP}\xspace}

\newcommand{\set}[1]{\ensuremath{\left\{#1\right\}}}

\newcommand{\floor}[1]{\ensuremath{\left\lfloor\!\!\begin{array}{l}#1\end{array}\!\!\right\rfloor}}
\newcommand{\paren}[1]{\ensuremath{\left(#1\right)}}
\newcommand{\prob}[1]{\ensuremath{P\!\left(#1\right)}}

\newcommand{\shorthand}[1]{{#1}\xspace}
\newcommand{\ie}{\shorthand{i.e.}}

\newcommand{\eg}{\shorthand{e.g.}}

\newcommand{\BigOh}[1]{\ensuremath{O\!\left(#1\right)}}

\begin{document}
\maketitle

\begin{abstract}
    We revisit the \emph{Strong Birthday Problem (\sbp)} introduced by DasGupta~\cite{DasGupta2005}, which asks for the minimum population $n$ required such that, with a probability of at least $1/2$, every individual in the group shares a birthday with at least one other person. Formally, we develop and analyze computational frameworks to determine the probability that in a group of $n$ people with birthdays distributed over $m$ days, each day either has two or more birthdays or is birthday-free. We derive both counting-based and probability-based recurrence relations to solve this problem and establish a novel connection to associated Stirling numbers of the second kind. This relationship is exploited to derive new, more efficient recurrences. Finally, we implement these recurrences using dynamic programming, provide analysis of their asymptotic complexities, and present numerical evaluations that demonstrate the practical efficiency and scalability of our proposed approaches.
\end{abstract}


\section{Introduction}\label{sec:introduction}
What is the probability that in a randomly chosen group of $n$ people, at least two will have the same birthday? This problem is famously known as the \emph{birthday problem} (see~\cite{Feller1968introductionV1}, page 33). While first published by Richard von Mises in 1939~\cite{vonMises1939}, it is generally attributed to Harold Davenport, who conceived of the problem approximately 12 years earlier but did not publish it at the time. 

A standard setup for the birthday problem assumes that birthdays are uniformly distributed over the days of the year. This includes several simplifying assumptions including but not limited to no selection bias (\eg, preference for a specific birthday or month), no interference between birth events (\eg, medical delays to accommodate hospital emergencies), and no multiple births such as twins. Under these conditions, only 23 people are required to achieve a probability greater than 1/2 that at least two individuals share a birthday. This result is often considered counterintuitive given that 366 people are required to reach a probability of 1 by the pigeonhole principle (disregarding leap years). A simple variant generalizes this to matching birthdays of couples, that is, wife-wife and husband-husband birthday matches; by replacing the 365-day year with $365^{2}$ possible birthday pairs, it can be shown that 431 couples are needed to reach a probability of 1/2 or more~\cite{DasGupta2005}. Another variant identifies the minimum number of people required for at least two to have birthdays within $d$ adjacent days with a probability of at least 1/2~\cite{Naus1968}. For $d = 1, 2, 3, 4, 5$, the corresponding minimum numbers of people are $n = 23, 14, 11, 9, 8$, respectively.

Several variants and extensions have been proposed to the birthday problem over the past eight decades~\cite{Naus1968,mosteller1962,Mathis1991,Nunnikhoven1992,DasGupta2005,inoue2008}. While many of them can be entertaining puzzles, their solutions, computationally efficient exact and approximation algorithms designed to compute those have been an active area of research. This stems from the fact that several real-world problems in the domains of computer science including hashing, cryptography~\cite{bellare2004}, communication networks and load balancing in distributed systems~\cite{eijkhout2019}, and problems in other domains including biology~\cite{Green2024} and medicine~\cite{Krawitz2015}, forensic science~\cite{Kaye2013Beyond}, statistics, social science and data science~\cite{gold2024}, can be modeled as one of the variants of the birthday problem. 

In this article, we consider the Strong Birthday Problem (\sbp) introduced in~\cite{DasGupta2005}. The problem asks to compute the minimum number of people we have to choose so that everyone has a shared birthday with probability at least 1/2. While closed-form formulas offer direct solutions, in this article, we explore dynamic programming (DP) recurrences as they provide a deeper understanding of the problem's combinatorial structure through explicit state transitions, generalize more effectively to variants with different constraints, and offer insight into the sequential evolution of probabilities as people are added. By establishing a novel connection to fundamental mathematical structures like partition theory through associated Stirling number of the second kind, we further derive efficient recurrences for \sbp. Our study systematically evaluates three distinct recurrence-based approaches against the direct formula, analyzing their theoretical complexity and practical performance to provide a comprehensive guide for algorithm selection.

The remainder of this article is organized as follows. Section~\ref{sec:strongbirthdayproblem} formally defines \sbp, followed by the derivation of a closed-form combinatorial formula in Section~\ref{subsec:combinatorialformula}. We then present three distinct algorithmic approaches: a counting-based recurrence relation in Section~\ref{subsec:countingrecurrence}, a direct probability-based recurrence in Section~\ref{subsec:directProbabilityRecurrence}, and a specialized recurrence exploiting the properties of associated Stirling numbers of the second kind in Section~\ref{subsec:assocStirlingNumbersAndConnToSBP}. Section~\ref{subsec:complexitySummary} provides a comparative analysis of the asymptotic time and space complexities for each derived method. We discuss our implementation strategy, utilizing both top-down memoization and bottom-up dynamic programming, and validate our results against the combinatorial baseline in Section~\ref{sec:impl}. Finally, Section~\ref{sec:extn} explores potential extensions and identifies related problems where our proposed framework may be generalized to broader classes of occupancy problems.

\section{The Strong Birthday Problem (\sbp)}\label{sec:strongbirthdayproblem}
\textbf{Problem Definition (Strong Birthday Problem):} Let $m \geq 1$ be the number of days in a year, and let $n \geq 1$ be the number of people whose birthdays are independently and uniformly distributed over the $m$ days. That is, the probability that the $i$th person's birthday is incident on the $j$th day of the year is $1/m$, for $1 \leq i \leq n$ and $1 \leq j \leq m$. For each day $d \in [m]$, let $b_d$ denote the number of people with birthdays incident on day $d$, where $\sum_{d=1}^{m} b_d = n$. 
Let $\gamma \in [0,1]$ be a probability threshold. Then the problem seeks to find the minimum number of people $n_{\min}$ such that with probability at least $\gamma$, every day has either zero birthdays or at least two birthdays, \ie, no day has exactly one birthday. Formally, we can write this as
\begin{eqnarray}\label{eqn:sbp}
n_{\min}(m, \gamma) &=& \min\left\{n \geq 1 : \Pr\left(\forall d \in [m]: b_d \neq 1\right) \geq \gamma\right\}.
\end{eqnarray}
In the above formulation, setting $\gamma = 1/2$, gives us the original statement of \sbp.

We approach the above problem by fixing $m, n$, and then finding the value of $\gamma$, and then doing a (binary) search on $n$ to find $n_{\min}$ so that $\gamma \geq 1/2$ (or any fixed probability threshold). Therefore, for the rest of the paper, without loss of generality, \emph{we use \sbp to refer to the problem of finding the probability $\gamma$ given the values of $m$ and $n$.}

\subsection{A Combinatorial Formula for \sbp}\label{subsec:combinatorialformula}
The paper~\cite{DasGupta2005} gave a formula to compute the probability for \sbp. However, we could not find a derivation of it. Therefore, as a warm-up exercise, we give a derivation of the combinatorial formula for \sbp. Specifically, we show that in a  population of $n$ people and a year consisting of $m$ days, the number of ways to form \emph{exactly} $k$ singleton birthdays---defined as days on which exactly one birthday occurs, is given by
\begin{eqnarray}\label{eqn:combinatorialFormula}
    N_{m, n, k} &=& \sum_{j = k}^{\min(m, n)} (-1)^{j - k} {j \choose k} {m \choose j} {n \choose j} j! (m - j)^{n - j}.
\end{eqnarray}

We give a proof of this in three steps: first establishing a formula for \emph{at least} $j$ singletons, then applying the inclusion-exclusion principle, and finally deriving the closed form formula.

\bigskip
\noindent\textbf{Step 1: Counting configurations with at least $j$ singleton birthdays.}

Let $M_{m,n,j}$ denote the number of ways to distribute $n$ people among $m$ days such that \emph{at least} $j$ specific days are singleton birthdays, \ie, each of these $j$ days has exactly one person.

To construct such a configuration:
\begin{enumerate}
    \item We choose $j$ specific birthdays from $m$ total days in the year in $m \choose j$ ways.
    \item We choose $j$ people from $n$ people whose birthdays will be incident on those selected days in $n \choose j$ ways.
    \item We establish a bijection between the $j$ chosen people and the $j$ chosen days. The number of such bijections is $j!$.
    \item Finally, we distribute the remaining $n - j$ people among the remaining $m - j$ available days. This can be performed in $(m - j)^{n - j}$ ways. Note that this process may produce additional singleton days; this is consistent with our objective to account for at least $j$ singletons.
\end{enumerate}

Therefore, the number of ways of forming at least $j$ singleton birthdays, \ie, the sum of configurations over all possible sets of $j$ specific days singletons, is given by
\begin{equation}
M_{m, n, j} = {m \choose j} {n \choose j} j! (m - j)^{n - j}.
\end{equation}

\noindent\textbf{Step 2: Applying the Inclusion-Exclusion principle.}

We now derive the number of configurations containing \emph{exactly} $k$ singleton birthdays, as opposed to \emph{at least $k$.} Let $S$ denote the set of all possible assignments of $n$ people to $m$ days. For any subset of days $D\subseteq [m]$ with cardinality $|D|=j$, let $C_{D}$ represent the set of configurations where every day in $D$ is a singleton birthday. While $|C_{D}|$ gives the number of configurations where at least the days in the specific set $D$ are singletons, the term $M_{m,n,j}$ accounts for all possible choices of $j$ days. Specifically, $M_{m,n,j}$ represents the sum of $|C_{D}|$ over all subsets $D$ of size $j$. Therefore,
\begin{equation}
M_{m,n,j} = \sum_{|D|=j} |C_D| = \binom{m}{j} \times (\text{ways to make } j \text{ specific days singletons}).
\end{equation}

By the inclusion-exclusion principle, the number of configurations with \emph{exactly} $k$ singleton birthdays, denoted by $S_k$, is obtained by starting with configurations that include $k$ specific singletons, subtracting those that have $k+1$ or more singletons, and adding back over-counted terms, and so on. This systematic addition and subtraction ensures that each configuration is counted precisely once. Formally,
\begin{equation}
|S_k| = \sum_{j = k}^{\min(m, n)} (-1)^{j - k} \binom{j}{k} \times (\text{configurations with at least $j$ singletons}),
\end{equation}
where the binomial coefficient $\binom{j}{k}$ accounts for the number of ways a configuration with $j$ singletons is overcounted among the $\binom{m}{k}$ possible sets of $k$ singletons, \ie, it counts the number of ways of choosing which $k$ of the $j$ \emph{at least singleton} days are the \emph{exactly $k$} singletons. The alternating sign $(-1)^{j-k}$ performs the inclusion-exclusion adjustment. Consequently, the total number of configurations with exactly $k$ singleton birthdays, $N_{m, n, k}$, is given by
\begin{equation}\label{eqn:inclusionExclusionFormula}
N_{m, n, k} = \sum_{j = k}^{\min(m, n)} (-1)^{j - k} {j \choose k} M_{m, n, j}.
\end{equation}

Note that the upper limit of the above summation is $\min(m,n)$ because the number of singleton birthdays cannot exceed the number of available days $m$ or the total number of people $n$.

\vspace{0.3cm}
\noindent\textbf{Step 3: Deriving the closed-form solution.}

Substituting the expression for $M_{m,n,j}$ from Equation~\ref{eqn:combinatorialFormula} into Equation~\ref{eqn:inclusionExclusionFormula}, we obtain the final closed-form solution for the number of ways of forming $k$ singleton birthdays:
\begin{eqnarray}
N_{m, n, k} &=& \sum_{j = k}^{\min(m, n)} (-1)^{j - k} {j \choose k} M_{m, n, j} \nonumber \\
            &=& \sum_{j = k}^{\min(m, n)} (-1)^{j - k} {j \choose k} {m \choose j} {n \choose j} j! (m - j)^{n - j}.\label{eqn:combinatorialFormulaForExactK}
\end{eqnarray}

This completes the derivation.

We note that for $0 \leq k \leq \min(n, m)$, the Equation~\ref{eqn:combinatorialFormulaForExactK} is valid, and for $k > m$ or $k > n$, $N_{m, n, k} = 0$ due to the impossibility of having more singleton birthdays than people or days.

\bigskip

\noindent\textbf{Computing the Probability:} The number of ways $n$ people can be assigned $m$ days as birthdays is $m^n$, which is the size of the sample space. Therefore, the probability of having $k$ singleton birthdays is
\begin{eqnarray}
    \prob{m, n, k} &=& \frac{1}{m^n} N_{m, n, k} = \frac{1}{m^n} \sum_{j = k}^{\min(m, n)} (-1)^{j - k} {j \choose k} {m \choose j} {n \choose j} j! (m - j)^{n - j}.\label{eqn:probabilityFormulaForExactK}
\end{eqnarray}

The above formula gives us the probability of having $k$ singleton birthdays, which itself can be interpreted as yet another variant of the birthday problem. Note that he remaining $m - k$ days can have 0, 2 or more than 2 (at most $n - k$) birthdays incident on them. 

Using $k = 0$ in the above formula, we obtain the probability that each person in the group of $n$ people shares her/his birthday with someone in the group:
\begin{eqnarray}
    \prob{m, n, 0} &=& \frac{1}{m^n} \sum_{j = 0}^{\min(m, n)} (-1)^{j} {m \choose j} {n \choose j} j! (m - j)^{n - j}.\label{eqn:probabilityFormulaForSBP}
\end{eqnarray}

\noindent\textbf{Complexity Analysis:} The formula in the Equation~\ref{eqn:probabilityFormulaForSBP} is computed by iterating the sum from $j = 0$ to $j = n$, and computing each term using the binomial coefficients, factorials and exponentials. Under the unit-cost arithmetic model, we have \BigOh{n} iterations with \BigOh{1} operations per iteration; therefore, the time complexity is \BigOh{n}. With arbitrary precision arithmetic, each operation has cost dependent on the bit-length of representation of the number; therefore, the time complexity is \BigOh{n \log^2(mn)}. The space complexity is \BigOh{1} under unit-cost arithmetic model, and it is \BigOh{\log^2(mn)} under bit-length model.

\subsection{A Recurrence for \sbp Based on Counting}\label{subsec:countingrecurrence}
Let $T(j, k, n, m)$ denote the number of ways in which $n$ birthdays are distributed over $m$ possible days in a year so that there are $j$ days with 2 or more birthdays incident on them and $k$ days with singleton birthdays. Therefore, $m - j - k$ days are birthday-free. Clearly, $n \geq 2 j + k$. We derive the recurrence by considering the birthday of the $n$th person. Given a valid configuration with $n - 1$ people, addition of $n$th person creates exactly one of the three mutually exclusive cases as given below.

\begin{enumerate}[{\bf {Case} 1:}]
    \item \textbf{Add to days with $\geq 2$ birthdays.} The $n$th person's birthday is incident on one of the $j$ days with 2 or more birthdays. Then, there are $T(j, k, n - 1, m)$ ways in which the first $n - 1$ birthdays are distributed so that there are $j$ days with 2 or more birthdays incident on them and $k$ days with singleton birthdays. Therefore, the number of ways this can be done is $j \cdot T(j, k, n - 1, m)$.

    \item \textbf{Add to singleton birthdays.} The $n$th person's birthday is incident on one of the $k + 1$ singleton birthdays. In this case, after $n$th person is included, the number of singletons birthdays goes from $k + 1$ to $k$ and the number of days with 2 or more birthday goes from $j - 1$ to $j$. Therefore, the number of ways this can be done is $(k + 1) \cdot T(j - 1, k + 1, n - 1, m)$.

    \item \textbf{Create a singleton birthday.} The $n$th person forms a singleton birthday. There are $T(j, k - 1, n - 1, m)$ ways in which the first $n - 1$ birthdays are distributed so that there $j$ days with 2 or more birthdays incident on them and $k - 1$ days with singleton birthdays. One of the remaining $m - j - (k - 1) = m - j - k + 1$ days can be used to realize the $n$th person's birthday. Therefore, the number of ways this can be done is $(m - j - k + 1) \cdot T(j, k - 1, n - 1, m)$.
\end{enumerate}

The above three cases partition all possible configurations giving the following recurrence:
\begin{eqnarray}
    T(j, k, n, m) &=& j \cdot T(j, k, n - 1, m) \nonumber \\
    && + (k + 1) \cdot T(j - 1, k + 1, n - 1, m) \nonumber \\
    && + (m - j - k + 1) \cdot T(j, k - 1, n - 1, m).
\end{eqnarray}

The base cases for the recurrence come from the following simple observations. When $n < 2j + k$, $T(j, k, n, m) = 0$ by the set up of the problem. When $j = k = n = 0$, $T(0, 0, 0, m) = 1$, for all $m \geq 0$, as there is one way of doing such an assignment -- by doing nothing. When $j < 0$ or $k < 0$, $T(j, k, n, m) = 0$. Also, when $m < j + k$, $T(j, k, n, m) = 0$.

Therefore, we can write the complete recurrence as
\begin{eqnarray}\label{eqn:countingRecurrence}
    T(j, k, n, m) &=& 
    \begin{cases}
        0, & \text{if $j < 0 \lor k < 0 \lor n < 2j + k$} \lor m < j + k; \\
        1, & \text{if $j = 0 \land k = 0 \land n = 0$}; \\
        j \cdot T(j, k, n - 1, m) \\ + (k + 1) \cdot T(j - 1, k + 1, n - 1, m) \\ + (m - j - k + 1) \cdot T(j, k - 1, n - 1, m), & \text{otherwise.} 
    \end{cases}
\end{eqnarray}

The above recurrence can be computed using dynamic programming.

For the strong birthday problem, $k = 0$, and $1 \leq j \leq \floor{\frac{n}{2}}$, and the number of ways the birthdays are formed so that each day have either 2 or more birthdays or is birthday-free is then
\begin{eqnarray}
    N &=& \sum_{j = 1}^{\floor{\frac{n}{2}}} T(j, 0, n, m),
\end{eqnarray}
and hence the probability of this happening is given by
\begin{eqnarray}
    \prob{m, n} &=& \frac{1}{m^n} N = \frac{1}{m^n} \sum_{j = 1}^{\floor{\frac{n}{2}}} T(j, 0, n, m).\label{eqn:sbpfirstprincrec}
\end{eqnarray}

\noindent \textbf{Computational Complexity:} The above formula can be computed using dynamic programming. Since in the recurrence $j \in [0, n / 2]$, $k \in [0, n]$ as it can grow during recursion due to the middle term in the recurrence, and $n$'s value ranges from 0 to $n$, the total number of states are \BigOh{n^3} where each state is computed in \BigOh{1} time under unit-cost model and assuming that memoization and retrieval time is \BigOh{1}. Therefore, the time complexity is \BigOh{n^3}. Bottom-up dynamic programming using rolling array that maintains only two layers of \BigOh{n^2} states, has the space complexity is \BigOh{n^2} and the time complexity \BigOh{n^3}.

\subsection{A Recurrence for \sbp Based on Direct Probability Computation}\label{subsec:directProbabilityRecurrence}
In Equation~\ref{eqn:sbpfirstprincrec}, we count the number of configurations first and then at the end perform the division by $m^n$. These operations needs arbitrary precision arithmetic that adds to the bit complexity. To mitigate this, we can derive a recurrence to directly compute and track the probabilities.

Let $P(j, k, n, m)$ denote the probability that $n$ birthdays distributed over $m$ days result in $j$ days with $\geq 2$ birthdays, $k$ singleton days, and $m - j - k$ birthday-free days. We derive the recurrence by considering the birthday of the $n$th person, given the valid probability states with $n - 1$ people. The probability transitions are from three different states with $n - 1$ people to the state with $n$ people, as given below.

\begin{enumerate}[{\bf {Case} 1:}]
    \item \textbf{Add to days with $\geq 2$ birthdays.} The $n$th person's birthday is incident on one of the $j$ days with $\geq 2$ birthdays. The probability of this happening is $j / m$ which contributes $(j / m) \cdot P(j, k, n - 1, m)$ to the total probability $P(j, k, n, m)$.

    \item \textbf{Add to singleton birthdays.} The $n$th person's birthday is incident on one of the $k + 1$ singleton birthdays. In this case, after $n$th person is included, the number of singletons birthdays goes from $k + 1$ to $k$ and the number of days with $\geq 2$ birthday goes from $j - 1$ to $j$. The probability that of this happening is $(k + 1) / m$ which contributes $(k + 1) / m \cdot P(j - 1, k + 1, n - 1, m)$ to the total probability $P(j, k, n, m)$.

    \item \textbf{Create a singleton birthday.} The $n$th person forms a singleton birthday by choosing one of the $m - j - (k - 1) = m - j - k + 1$ available birthday-free days. The probability of this happening is $(m - j - k + 1) / m$ which contributes $(m - j - k + 1) / m \cdot P(j, k - 1, n - 1, m)$ to the total probability $P(j, k, n, m)$.
\end{enumerate}

Summing these three mutually exclusive cases, and adding the base cases as we did in deriving Equation~\ref{eqn:countingRecurrence} yields the following recurrence:
\begin{eqnarray}
    P(j, k, n, m) &=& 
    \begin{cases}
        0, & \text{if $j < 0 \lor k < 0 \lor n < 2j + k$} \lor m < j + k; \\
        1, & \text{if $j = 0 \land k = 0 \land n = 0$}; \\
        \frac{j}{m} \cdot P(j, k, n - 1, m) \\ + \frac{k + 1}{m} \cdot P(j - 1, k + 1, n - 1, m) \\ + \frac{m - j - k + 1}{m} \cdot P(j, k - 1, n - 1, m), & \text{otherwise.}\label{eqn:probabilityRecurrence}
    \end{cases}
\end{eqnarray}

For \sbp, the probability is computed using the following equation:
\begin{eqnarray}
    P(m, n) &=& \sum_{j = 1}^{\floor{\frac{n}{2}}} P(j, 0, n, m).\label{eqn:sbpProbRecurrence}
\end{eqnarray}

We note that the two recurrences in Equation~\ref{eqn:countingRecurrence} and Equation~\ref{eqn:probabilityRecurrence} are related by
\begin{eqnarray}
    P(j, k, n, m) &=& \frac{T(j, k, n, m)}{m^n},
\end{eqnarray}
which can be proved straightforwardly by using induction (not shown here).

\noindent \textbf{Computational Complexity:} The bottom-up dynamic programming implementation with rolling array (only two layers of \BigOh{n^2} states to store the current and previous states are needed) has the time and space complexities \BigOh{n^3} and \BigOh{n^2}, respectively. However, since the probability recurrence uses standard floating-point arithmetic as opposed to arbitrary precision arithmetic, the arithmetic operations are done in \BigOh{1} time which benefits from hardware-accelerated operations, giving significant speedup in practical implementations.

\subsection{Associated Stirling Numbers of the Second Kind and It's Connection to the Strong Birthday Problem}\label{subsec:assocStirlingNumbersAndConnToSBP}
We show that associated Stirling numbers of the second kind can be used to derive novel recurrences for \sbp that are both time and space efficient. Specifically, we show that the number of unique states are \BigOh{n^2}, and only \BigOh{n} states need to be stored in a bottom-up rolling array implementation, giving us \BigOh{n^2} time and \BigOh{n} space comepexities. Since the recurrence is independent of $m$ during the dynamic programming computation, it enables precomputation of Stirling numbers. The recurrence we derive here has cleaner state space than the method based on counting and direct probability computation, and has better asymptotic complexity bounds.

\subsubsection{Associated Stirling Numbers of the Second Kind}\label{subsubsec:assocstr}
Stirling numbers of the \emph{second kind}, denoted by ${n \brace k}$, counts the number of ways to partition $n$ distinct objects \set{a_1, \ldots, a_n} into $k$ nonempty partitions (subsets). For example, ${4 \brace 2} = 7.$

We can write a recurrence for ${n \brace k}$ as follows (see~\cite{graham94} page 259):

\begin{eqnarray}\label{eqn:strilingSecondKind}
    {n \brace k} &=& 
    \begin{cases}
        0, & \text{if $(n = 0 \land k > 0) \lor (n > 0 \land k = 0)$}; \\
        1, & \text{if $(n = 0 \land k = 0) \lor (n > 0 \land k = 1)$}; \\
        {n - 1 \brace k - 1} + k {n - 1 \brace k}, & \text{otherwise.}
    \end{cases}
\end{eqnarray}

The \emph{$r$-associated} Stirling number of the second kind, denote by ${n \brace k}_{\geq r}$ is the number of ways of partitioning $n$ distinct objects, \set{a_1, \ldots, a_n}, into $1 \leq k \leq \floor{\frac{n}{r}}$ \emph{unlabeled} partitions where each partition has size at least $r \geq 2$. Note that this is a restriction over the Stirling number of the second kind where $r \geq 1$. For example, ${7 \brace 2}_{\geq 3} = 35$ and ${7 \brace 2}_{\geq 2} = 105$.

We can derive a recurrence for ${n \brace k}$ using the following insight. The $n$th object $a_n$ can be considered in two mutually exclusive ways as follows.

\begin{enumerate}[{\bf {Case} 1:}]
    \item \textbf{Forms a new partition.} The $n$th object's insertion creates a new partition of size exactly $r$. The number of ways in which $r - 1$ objects can be chosen from $n - 1$ objects is $n - 1 \choose r - 1$, to which the $n$th object is added. The remaining $n - r$ objects can be partitioned into $k - 1$ partitions of size at least $r$ in ${n - r \brace k - 1}_{\geq r}$ ways. Therefore, the number of ways are ${n - 1 \choose r - 1} {n - r \brace k - 1}_{\geq r}$.

    \item \textbf{Does not form a new partition.} The $n$th object joins an existing partition. Then, there are ${n - 1 \brace k}_{\geq r}$ ways of partitioning $n - 1$ objects into $k$ partitions of size at least $r$. The $n$th objects can go into any of these $k$ partitions. Therefore, the number of ways are $k {n - 1 \brace k}_{\geq r}$.
\end{enumerate}

Therefore, the recurrence is given by
\begin{eqnarray}
    {n \brace k}_{\geq r} &=& k {n - 1 \brace k}_{\geq r} + {n - 1 \choose r - 1} {n - r \brace k - 1}_{\geq r}.
\end{eqnarray}

The base cases for the recurrence come from the following simple observations: ${0 \brace 0}_{\geq r} = 1$, and ${n \brace k}_{\geq r} = 0$ when $n \leq 0$ or $k \leq 0$ or $n < kr$. Also, ${n \brace n}_{\geq r} = 1$. 

Therefore, we can write the complete recurrence as:
\begin{eqnarray}\label{eqn:assocStirlingSecondKind}
    {n \brace k}_{\geq r} &=& 
    \begin{cases}
        0, & \text{if $(n \leq 0 \land k > 0) \lor (n > 0 \land k \leq 0) \lor (n < kr)$}; \\
        1, & \text{if $n = 0 \land k = 0$}; \\
        k {n - 1 \brace k}_{\geq r} + {n - 1 \choose r - 1} {n - r \brace k - 1}_{\geq r}, & \text{otherwise.}
    \end{cases}
\end{eqnarray}

For $r = 1$, the above recurrence reduces to Stirling numbers of the second kind (see Equation~\ref{eqn:strilingSecondKind}). For $r = 2$, which is relevant to \sbp, the above recurrence reduces to
\begin{eqnarray}\label{eqn:assocStirlingSecondKindForRvalue2}
    {n \brace k}_{\geq 2} &=& k {n - 1 \brace k}_{\geq 2} + (n - 1) {n - 2 \brace k - 1}_{\geq 2}.
\end{eqnarray}

See~\cite{Comtet1974,Howard1980,Zhao2008,Komatsu2015} for more on the associated Stirling numbers of the second kind and related identities. We will use the above-derived recurrence for the associated Stirling numbers of the second kind to formulate a recurrence for \sbp.

\subsubsection{A Recurrence for \sbp Using the Associated Stirling Numbers of the Second Kind}\label{subsubsec:recurrenceusingassocstr}
When the partitions are \emph{labeled,} as in case of the birthday problem where days are distinguishable, the number of ways of putting $n$ labeled objects into $k$ labeled partitions such that each partition has at least $r$ objects can be obtained by considering all $k!$ permutations of the $k$ partition labels. Therefore,
\begin{eqnarray}
    T_r(n, k) &=& k! {n \brace k}_{\geq r} \nonumber \\
        &=& \begin{cases}
            0, & \text{if $(n \leq 0 \land k > 0) \lor (n > 0 \land k \leq 0) \lor (n < kr)$}; \\
            1, & \text{if $n = 0 \land k = 0$}; \\
            k! \paren{k {n - 1 \brace k}_{\geq r} + {n - 1 \choose r - 1} {n - r \brace k - 1}_{\geq r}}, & \text{otherwise.}
        \end{cases}
\end{eqnarray}

This recurrence has a nice interpretation as one more variant of the birthday problem: \emph{How many ways are there so that $n$ people's birthdays can be assigned to $k$ days so that each day has $r \geq 2$ person's birthdays assigned to it?} The corresponding probability problem can be stated accordingly.

For $r = 2$, the above recurrence reduces to
\begin{eqnarray}\label{eqn:sbpAssocCountWaysRvalue2}
    T_2(n, k) &=& k! {n \brace k}_{\geq 2} = k! \paren{k {n - 1 \brace k}_{\geq 2} + (n - 1) {n - 2 \brace k - 1}_{\geq 2}}.
\end{eqnarray}

For the strong birthday problem, we can derive a recurrence using the above recurrences in Equation~\ref{eqn:sbpAssocCountWaysRvalue2}, by counting no singleton birthdays, \ie, the number of favorable outcomes, using the following two steps.

\begin{enumerate}[{\bf {Step} 1:}]
    \item \textbf{Choose $k$ specific days from $m$ days.} The number of ways of choosing $k$ specific days from $m$ days is $\binom{m}{k}$, where $1 \leq k \leq \floor{\frac{n}{2}}$.
    
    \item \textbf{Assign $n$ People to $k$ labeled days.} The number of ways of assigning $n$ people's birthdays to $k$ specific days such that each day has $\geq 2$ person's birthdays incident on it is $T_2(n, k) = k! {n \brace k}_{\geq 2}$.
\end{enumerate}

Therefore, the total number of ways is given by
\begin{eqnarray}\label{eqn:sbpFinalAssocCountWaysRvalue2}
    N &=& \sum_{k = 1}^{\floor{\frac{n}{2}}} \binom{m}{k} T_2(n, k) = \sum_{k = 1}^{\floor{\frac{n}{2}}} \binom{m}{k} k! {n \brace k}_{\geq 2}.
\end{eqnarray}

Since the total number of ways in which $n$ people can be assigned $m$ days as birthdays is $m^n$, which is the size of the sample space, the probability that each person shares her/his birthday with someone in the group of $n$ people is
\begin{eqnarray}\label{eqn:sbpProbabilityAssocSterlingRvalue2}
    \prob{m, n, r = 2} &=& \frac{1}{m^n} N = \frac{1}{m^n} \sum_{k = 1}^{\floor{\frac{n}{2}}} {m \choose k} k! {n \brace k}_{\geq 2}.
\end{eqnarray}

We can have a simple generalization to the strong birthday problem when we ask the following question: \emph{What is the probability that in a group of $n$ people and $m$ days in a year, each day has at least $r \geq 2$ birthdays incident on it or is birthday-free?} In other words, what is the probability that each person's birthday is shared by at least $r - 1$ other people for $r \geq 2$? The solution can be written straightforwardly, by observing that $1 \leq k \leq \floor{\frac{n}{r}}$, and is given by
\begin{eqnarray}\label{eqn:sbpProbabilityGeneralizedFormula}
    \prob{m, n, r} &=& \frac{1}{m^n} \sum_{k = 1}^{\floor{\frac{n}{r}}} {m \choose k} k! {n \brace k}_{\geq r}.
\end{eqnarray}

The recurrence in Equation~\ref{eqn:sbpProbabilityAssocSterlingRvalue2} can be computed using dynamic programming, and has several significant advantages over the previous recurrences we derived (see Equation~\ref{eqn:countingRecurrence} and Equation\ref{eqn:probabilityRecurrence}). The state space described by ${n \brace k}_{\geq 2}$ only depends on $n$ and is independent of $m$. The recurrence in Equation~\ref{eqn:assocStirlingSecondKindForRvalue2} comprises of two terms instead of three as in the previous recurrences and have small integer multipliers. Most importantly, since ${n \brace k}_{\geq 2}$ computation is independent of $m$, they can be precomputed for given $n$ and can be reused across different values of $m$ which only appears in the final summation phase not during the dynamic programming recursion. Finally, since the division by $m^n$ and multiplication by $\binom{m}{k}$ occurs in the final summation step, and the Stirling number recursion operates on values whose growth depends only on $n$ and $k$ with predictable intermediate sizes of large integers there by avoiding manipulations of enormous integers required to implement Equation~\ref{eqn:countingRecurrence}, and also avoids the repeated division by $m$ required by to implement Equation\ref{eqn:probabilityRecurrence}.

\noindent \textbf{Computational Complexity:} The dynamic programming implementation using top-down memoization has \BigOh{n^2} states where each state requires \BigOh{1} operations under unit-cost arithmetic model. Therefore, both time and space complexities are \BigOh{n^2}. The bottom-up dynamic programming implementation using (auxiliary) rolling arrays has time complexity \BigOh{n^2}, and space complexity \BigOh{n} since only two rows of \BigOh{n} states need to be stored to compute the current row.

\subsection{Summary of Complexity of Algorithms for \sbp}\label{subsec:complexitySummary}
Table~\ref{tab:complexity} summarizes the asymptotic time and space complexities of the various methods derived for solving \sbp. Compared to the implementations of the recurrences in Equation~\ref{eqn:sbpfirstprincrec} and Equation~\ref{eqn:sbpProbRecurrence}, which require \BigOh{n^3} time, the recurrence in Equation~\ref{eqn:sbpProbabilityAssocSterlingRvalue2} is asymptotically faster with a time complexity of \BigOh{n^2}. This efficiency gain extends to memory usage as well. The top-down memoization of the former recurrences requires \BigOh{n^3} space, which can be reduced to \BigOh{n^2} using bottom-up dynamic programming with rolling arrays. In contrast, Equation~\ref{eqn:sbpProbabilityAssocSterlingRvalue2} requires only \BigOh{n^2} space with top-down memoization, which further reduces to \BigOh{n} space when implemented via bottom-up dynamic programming using rolling arrays.

\begin{table}[h]
\centering
\caption{Complexity Comparison of Algorithms for \sbp}
\label{tab:complexity}
\begin{tabular*}{0.8\textwidth}{@{\extracolsep{\fill}}lcc@{}}
\toprule
\textbf{Algorithm} & \textbf{Time Complexity} & \textbf{Space Complexity} \\
\midrule
\textbf{Combinatorial Formula (Equation~\ref{eqn:probabilityFormulaForSBP})} & $\BigOh{n}$ & $\BigOh{1}$ \\
\midrule
\multicolumn{3}{@{}l}{\textbf{Counting Recurrence (Equation~\ref{eqn:sbpfirstprincrec})}} \\
\quad Memoization & $\BigOh{n^3}$ & $\BigOh{n^3}$ \\
\quad Bottom-Up DP & $\BigOh{n^3}$ & $\BigOh{n^2}$ \\
\midrule
\multicolumn{3}{@{}l}{\textbf{Probability Recurrence (Equation~\ref{eqn:sbpProbRecurrence})}} \\
\quad Memoization (float) & $\BigOh{n^3}$ & $\BigOh{n^3}$ \\
\quad Bottom-Up (float) & $\BigOh{n^3}$ & $\BigOh{n^2}$ \\
\midrule
\multicolumn{3}{@{}l}{\textbf{Stirling Recurrence (Equation~\ref{eqn:sbpProbabilityAssocSterlingRvalue2})}} \\
\quad Memoization & $\BigOh{n^2}$ & $\BigOh{n^2}$ \\
\quad Bottom-Up DP & $\BigOh{n^2}$ & $\BigOh{n}$ \\
\bottomrule
\end{tabular*}
\end{table}

\section{Results and Discussion}\label{sec:impl}
The different algorithms for computing \sbp were implemented in Python 3.x. To facilitate rigorous verification, the combinatorial formula, the counting-based recurrence, and the associated Stirling number recurrence were implemented using the \texttt{mpmath} library with 1000 decimal places of precision (\texttt{mp.dps = 1000}), producing identical results across all test cases. In contrast, our implementation of the probability-based recurrence utilizes standard 64-bit IEEE 754 floating-point arithmetic. This approach demonstrates a significant performance gains and exhibits lower sensitivity to rounding errors compared to counting-based methods, making it ideal for scenarios where standard precision suffices. All experiments were conducted on a macOS-based system.

Table~\ref{tab:algorithmMasterComparison} summarizes the asymptotic performances and numerical characteristics of the combinatorial formula and the optimized bottom-up dynamic programming algorithms. While the combinatorial formula offers the lowest theoretical time complexity, its reliance on arbitrary-precision libraries for stability makes the probability-based recurrence more practical for simulations using native floating-point hardware. Stirling recurrence (Equation~\ref{eqn:sbpProbabilityAssocSterlingRvalue2}) offers best of both the worlds with arbitrary precision. Although it uses \texttt{mpmath} arithmetic functions, due to its cleaner state space than the other two methods, and better asymptotic time and space complexities compared to the counting and probability recurrences, it provides the best balance when high precision is needed. Furthermore, Stirling recurrence reveals connections to partition theory through associated Sterling numbers of the second kind which is exploited to derive the recurrence.

\begin{table}[h]
\centering
\caption{Comparison of Algorithms for \sbp: Complexity, Precision, and Stability}
\label{tab:algorithmMasterComparison}
\begin{tabularx}{0.8\textwidth}{@{}X l c c l@{}}
\toprule
\textbf{Algorithm} & \textbf{Precision} & \textbf{Time} & \textbf{Space\textsuperscript{$\dagger$}} & \textbf{Stability} \\ 
\midrule
\textbf{Combinatorial Formula (Equation~\ref{eqn:probabilityFormulaForSBP})} & High (\texttt{mpmath}) & $\BigOh{n}$ & $\BigOh{1}$ & Low\textsuperscript{*} \\
\addlinespace
\textbf{Counting Recurrence (Equation~\ref{eqn:sbpfirstprincrec})} & High (\texttt{mpmath}) & $\BigOh{n^3}$ & $\BigOh{n^2}$ & Moderate \\
\addlinespace
\textbf{Probability Recurrence (Equation~\ref{eqn:sbpProbRecurrence})} & Standard & $\BigOh{n^3}$ & $\BigOh{n^2}$ & High \\
\addlinespace
\textbf{Stirling Recurrence (Equation~\ref{eqn:sbpProbabilityAssocSterlingRvalue2})} & High (\texttt{mpmath}) & $\BigOh{n^2}$ & $\BigOh{n}$ & Moderate \\
\bottomrule
\addlinespace
\multicolumn{5}{l}{\footnotesize \textsuperscript{$\dagger$} Space complexity reflects optimized bottom-up dynamic programming implementations.} \\
\multicolumn{5}{l}{\footnotesize \textsuperscript{*} Requires arbitrary precision (\texttt{mpmath}) due to large intermediate values like $m^n$ and $n!$.} \\
\end{tabularx}
\end{table}

We selected 6 test cases to examine two critical probability thresholds across different scales of $m$:
\begin{itemize}
    \item Cases approaching $\approx 50\%$ probability: $(m,n) = (10, 41), (50, 304), (100, 690)$; and
    \item Cases approaching $\approx 99.9\%$ probability: $(m,n) = (10, 112), (50, 665), (100, 1410)$.
\end{itemize}

Table~\ref{tab:sbpAlgoImplementationComparisons} provides a comprehensive runtime comparison of the four implemented methods across varying scales of $m$ and $n$. The combinatorial formula based implementation with time complexity of $\BigOh{n}$ is much faster than the dynamic programming methods as expected, although numerical stability can be an issue for large values of $m$ and $n$. These results highlight the dramatic performance gains of closed-form solutions. By bypassing the computational depth of recurrences, these solutions offer both faster execution and simpler implementation paths.

\begin{table}[h]
\centering
\caption{Performance Comparison of Algorithms for \sbp: Runtime (Seconds)}
\label{tab:sbpAlgoImplementationComparisons}
\sisetup{table-format=2.3} 
\newcolumntype{Y}{>{\centering\arraybackslash}X} 

\begin{tabularx}{\textwidth}{@{} c cc SSS S @{}}
\toprule
 & \multicolumn{2}{c}{\textbf{Parameters}} & \multicolumn{4}{c}{\textbf{Runtime (s)}} \\
\cmidrule(lr){2-3} \cmidrule(l){4-7}
\textbf{Test} & $m$ & $n$ & {\textbf{Combinatorial Formula}} & {\textbf{Counting Rec.}} & {\textbf{Probability Rec.}} & {\textbf{Stirling Rec.}} \\
 & & & {\footnotesize\textbf{(Equation~\ref{eqn:probabilityFormulaForSBP})}} & {\footnotesize\textbf{(Equation~\ref{eqn:sbpfirstprincrec})}} & {\footnotesize\textbf{(Equation~\ref{eqn:sbpProbRecurrence})}} & {\footnotesize\textbf{(Equation~\ref{eqn:sbpProbabilityAssocSterlingRvalue2})}} \\
 & & & {\footnotesize\textbf{(mpmath)}} & {\footnotesize\textbf{(mpmath)}} & {\footnotesize\textbf{(float)}} & {\footnotesize\textbf{(mpmath)}} \\
\midrule
1 & 10 & 41 & 0.004 & 0.019 & 0.002 & 0.004 \\
2 & 50 & 304 & 0.027 & 4.047 & 0.350 & 0.165 \\
3 & 100 & 690 & 0.066 & 42.073 & 3.242 & 0.885 \\
4 & 10 & 112 & 0.008 & 0.061 & 0.007 & 0.022 \\
5 & 50 & 665 & 0.061 & 9.987 & 1.034 & 0.807 \\
6 & 100 & 1410 & 0.176 & 91.624 & 8.583 & 4.040 \\
\midrule
\multicolumn{3}{@{}l}{\textbf{Complexity}} & {$\BigOh{n}$} & {$\BigOh{n^3}$} & {$\BigOh{n^3}$} & {$\BigOh{n^2}$} \\
\multicolumn{3}{@{}l}{\textbf{Performance Ratio\textsuperscript{*}}} & {$1\times$} & {$520.6\times$} & {$48.8\times$} & {$23.0\times$} \\
\bottomrule
\addlinespace
\multicolumn{7}{l}{\footnotesize \textsuperscript{*} Performance Ratio is calculated relative to the Combinatorial Formula for Test 6.}
\end{tabularx}
\end{table}

All three recurrence relations, while computationally slower than the combinatorial formula, remain valuable for theoretical understanding. For example, the optimal substructures and the overlapping subproblems are not observed through the combinatorial formula; therefore, a limited insight into the structure of \sbp is obtained from the combinatorial formula. In contrast, the three recurrences derived here show the different combinatorial recursive structures of the same problem, respectively, through counting, probability flow and establishing a profound connection to partition theory. Therefore, we use the implementation of the combinatorial formula primarily as a baseline for correctness verification.

For the test case ($m = 100$ and $n = 1410$), we give the speed up ratios with respect to the run-time of the combinatorial formula. Among the recurrence-based methods, the Stirling recurrence shows the most significant advantage; for $n=1410$ (Test 6), its $\BigOh{n^2}$ complexity reduces runtime to $4.04$~seconds, yielding a $22.7$-fold speedup over the $\BigOh{n^3}$ counting recurrence while maintaining high precision via \texttt{mpmath}. Additionally, the probability recurrence remains competitive for smaller $n$ by leveraging 64-bit floating-point arithmetic, demonstrating the trade-off between computational throughput and numerical precision.

Table~\ref{tab:strongBirthday} presents the minimum population size $n$ required to ensure that every individual shares a birthday with at least one other person with probabilities $P \geq 0.5$ and $P \geq 0.999$. Our various implementations yielded identical results across all test cases, verifying the correctness of the derived approaches. For 
$m > 200$, the high time and space complexities of the counting and probability recurrences became computationally prohibitive; consequently, results for these larger scales were computed and cross-verified exclusively using the combinatorial formula (Equation~\ref{eqn:probabilityFormulaForSBP}) and the Stirling recurrence (Equation~\ref{eqn:sbpProbabilityAssocSterlingRvalue2}). 

Notably, while the classical birthday problem requires only $n \approx \sqrt{m}$ individuals to find a single pair, the Strong Birthday Problem requires a significantly larger population, scaling approximately as $n \sim m \log m$, \ie, $n = \BigOh{m \log m}$, to ensure that no individual remains a singleton. This asymptotic behavior is similar to that of the Coupon Collector's Problem~\cite{ErdosRenyi1961, Motwani1995} since the occupancy distribution~\cite{Feller1968introductionV1} in both the problems uses Stirling number of the second kind (and their restrictions), which forms the basis for our recurrences.

\begin{table}[h]
\centering
\caption{Minimum number of people $n$ required for a given $m$ such that every individual shares a birthday with at least one other person in the group with probability $P \geq 0.5$ and $P \geq 0.999$}
\label{tab:strongBirthday}
\begin{tabularx}{0.85\textwidth}{@{} X c c @{}}
\toprule
\textbf{Days ($m$)} & \textbf{$n$ at $P \geq 0.5$} & \textbf{$n$ at $P \geq 0.999$} \\
\midrule
10   & 41   & 112   \\
50   & 304  & 665   \\
100  & 690  & 1410  \\
200  & 1541 & 2975  \\
364  & 3054 & 5653  \\
365  & 3064 & 5669  \\
366  & 3073 & 5686  \\
400  & 3399 & 6253  \\
500  & 4375 & 7937  \\
1000 & 9528 & 16619 \\
\bottomrule
\end{tabularx}
\end{table}

\section{Significance and Extensions}\label{sec:extn}
Beyond its theoretical appeal, \sbp has high-impact applications in diverse fields. As discussed in \cite{DasGupta2005}, the problem provides a framework for fingerprint matching in criminology and sociological modeling. In computer science, \sbp is related to 
$k$-anonymity and honey encryption, where the elimination of singletons is required to ensure that every record in a database is indistinguishable from a larger anonymity set. In distributed systems and load balancing, \sbp models the point at which no server node is left with a lone task, ensuring redundancy and mirroring across all occupied resources. Furthermore, the problem finds relevance in job scheduling and retail operations; for example, it can model customer queues in a shopping store where management seeks to estimate the population $n$ required such that each active queue has a size $r \geq 2$. Such estimations are vital for making real-time decisions regarding merging or opening new service channels.

Several natural extensions of \sbp remain fertile ground for future research. A primary direction is the development of fast approximation algorithms to compute $n$ for arbitrary  $m$ and probability thresholds, particularly for scales where \BigOh{n^2} recurrences become computationally expensive. Additionally, relaxing the uniform probability assumption to account for skews or cases where individuals have non-identical probability distributions across the year would provide a more robust model for real-world phenomena.

\section*{Acknowledgments} I am thankful to Dr.~Chittaranjan Tripathy for introducing me to discrete mathematics, birthday problem and then suggesting me that \sbp can be an interesting problem to work on. I am thankful to Dr.~Shashidhara K.~Ganjugunte and Dr.~Pouyan Shirzadian for numerous valuable comments on the manuscript.
\bibliographystyle{unsrt}
\bibliography{references}  

\end{document}